\documentclass[12pt]{article}
\usepackage{amscd,amssymb,epic,amsmath}

\usepackage[dvips]{graphicx}
\usepackage[mathscr]{eucal}
\begin{document}
%\sloppy \raggedbottom \setcounter{page}{1}
%%%%%%%%%%%%%%%%%%%%%%%%%%

%\setcounter{figure}{0} \setcounter{equation}{0} \setcounter{footnote}{0}
%\setcounter{table}{0} \setcounter{section}{0}


% use the thanks command within \title, \author or \address for footnotes:
% \title{Title} or  \title{Title\thanks{...}}

\title{On some quantum bounded symmetric domains}\author{Olga Bershtein}

%\address{Institute for Low Temperature Physics and Engineering, Lenin ave.
%47, Kharkov 61103, Ukraine, bershtein@ilt.kharkov.ua}

\maketitle

\begin{abstract}
In the framework of quantum group theory we obtain a noncommutative analog
for the algebra of functions in a bounded symmetric domain, endowed with a
whole symmetry. Also we provide a construction for its faithfull irreducible
representation and an invariant integral over the bounded symmetric domain.
\end{abstract}

\newtheorem{proposition}{Proposition}
\newtheorem{theorem}{Theorem}

\section{Introduction}

Recall some well-known facts on bounded symmetric domains. We focus on a
series of bounded symmetric domains $\mathbb D$ from the well-known Cartan
list.

%Cartan obtained the classification of irreducible bounded symmetric domains
%(up to isomorphisms) given by Dynkin diagrams with a marked vertex.
%Harish-Chandra proved that every irreducible bounded symmetric domain
%$\mathbb D$ embeds into a finite dimensional normed vector space as the unit
%ball. Consider a case from the Cartan list which corresponds to the $C_n$-type
%Dynkin diagram with $n$-th marked vertex.

Let $$\mathfrak{a}=\begin{pmatrix}
2 & -1 & 0 & 0 & \ldots & 0 & 0 & 0\\
-1 & 2 & -1 & 0 & \ldots & 0 & 0 & 0\\
\ldots \\
0 & 0 & 0 & 0 & \ldots & -1 & 2 & -2 \\
0 & 0 & 0 & 0 & \ldots & 0 & -1 & 2 \\
\end{pmatrix}.$$

The Lie algebra $\mathfrak{g}=\mathfrak{sp}_{2n}(\mathbb C)$ is isomorphic
to the Lie algebra with generators $e_i,f_i,h_i,i=1,...,n$ and relations
$$
[h_i,h_j]=0, \quad [h_i,e_j]=a_{ij}e_j,\quad [h_i,f_j]=-a_{ij}f_j,$$
$$[e_i,f_j]=\delta_{ij}h_i \quad i,j=1,...,n,
$$
together with Serre's relations. The linear span $\mathfrak{h}$ of
$h_1,h_2,\ldots,h_n$ is a Cartan subalgebra, and the linear functionals
$\alpha_1,\alpha_2,\ldots,\alpha_n$ on $\mathfrak{h}$ defined by
$$\alpha_j(h_i)=a_{ij}, \qquad i,j=1,2,\ldots,n$$ form a system of
simple roots for the Lie algebra $\mathfrak{sp}_{2n}(\mathbb C)$.

Let $h_0 \in\mathfrak{h}$ be the element given by $$\alpha_{n}(h_0)=2,
\qquad \alpha_j(h_0)=0,\; j < n.$$ It is easy to prove that
$h_0=h_1+2h_2+...+nh_n$.

Let $\mathfrak{k} \subset \mathfrak{g}$ be the Lie subalgebra generated by
$$e_i,\;f_i,\quad i\neq n, \qquad \quad h_i,\; i=1,...,n.$$

The pair $(\mathfrak{g},\mathfrak{k})$ is Hermitian symmetric, i.e.
$\mathfrak{g}$ is equipped with a grading
$\mathfrak{g}=\mathfrak{p}^-\oplus\mathfrak{k}\oplus\mathfrak{p}^+$, where
\begin{align*}\label{par_type}
&\mathfrak{p}^{\pm}=\{\xi\in\mathfrak{g}|\:[h_0, \xi]=\pm 2\xi \},
\\& \mathfrak{k}=\{\xi\in\mathfrak{g}|\:[h_0,\xi]=0\}.
\end{align*}

Note that $\mathfrak{p}^-$ is isomorphic to the normed vector space of
symmetric complex $n \times n$-matrices with the operator norm.
Harish-Chandra proved that an irreducible bounded symmetric domain $\mathbb
D$ can be embedded into the normed vector space $\mathfrak{p}^-$ as the unit
ball.

In this paper we consider quantum analogs for the algebra $\mathbb
C[\mathfrak{p}^-]$ of holomorphic polynomials on $\mathfrak{p}^-$ and the
algebra $\mathrm{Pol}(\mathfrak{p}^-)$ of polynomials on
$\mathfrak{p}^-_{\mathbb R}$ and give some results on representation theory
on $\mathbb D$.

\section{Algebras $\mathbb C[\mathfrak{p}^-]_q$ and
$\mathrm{Pol}(\mathfrak{p}^-)_q$}
 In the sequel $q \in (0,1)$,
$\mathbb C$ is the ground field, and all the algebras are assumed
associative and unital.

Let $d_i, i=1,...,n$ be coprime numbers that symmetrize the Cartan matrix
$\mathfrak{a}$. One can check that $d_i=1$ for $i=1,...,n-1$ and $d_n=2.$

Denote by $U_q \mathfrak{g}=U_q \mathfrak{sp}_{2n}$ a Hopf algebra with
generators $K_i$, $K_i^{-1}$, $E_i$, $F_i$, $i=1,\ldots,n$, and relations
$$K_iK_j=K_jK_i,\quad K_iK_i^{-1}=K_i^{-1}K_i=1,$$
\begin{equation*}
K_iE_j=q^{d_ia_{ij}}E_jK_i,\quad K_iF_j=q^{-d_ia_{ij}}F_jK_i, \quad
i,j=1,...,n
\end{equation*}
$$
E_iF_j-F_jE_i=\delta_{ij}\,\frac{K_i-K_i^{-1}}{q^{d_i}-q^{-d_i}},
$$
together with $q$-analogs of the well-known Serre relations \cite{Klimyk}.

The coproduct $\Delta$, the counit $\varepsilon$ and the antipod $S$ are
defined as follows:
\begin{equation*}\label{comult_D}
  \Delta(E_i)=E_i\otimes 1+K_i\otimes E_i,\, \,
  \Delta(F_i)=F_i\otimes K_i^{-1}+1\otimes F_i, \, \,
  \Delta(K_i)=K_i\otimes K_i,
\end{equation*}
\begin{equation*}\label{antipode_DJ}
  S(E_i)=-K_i^{-1}E_i,\quad
  S(F_i)=-F_iK_i,\quad
  S(K_i)=K_i^{-1},
\end{equation*}
$$
\varepsilon(E_i)=\varepsilon(F_i)=0,\quad \varepsilon(K_i)=1.
$$

Denote by $U_q \mathfrak{k} \subset U_q \mathfrak{g}$ a Hopf subalgebra
generated by
$$
E_j,\,F_j, \quad j<n \qquad \text{ and } \qquad K^{\pm 1}_i, \quad
i=1,...,n.
$$

$U_q \mathfrak{g}$ can be equipped with the involution $*$ given by
$$
(K_j^{\pm})^*=K^{\pm}_j, \quad j=1,...,n.
$$
$$
E_j^*=
\begin{cases}
K_jF_j, & j< n,
\\ -K_jF_j, & j=n,
\end{cases} \quad
F_j^*=
\begin{cases}
E_jK_j^{-1}, & j< n,
\\ -E_jK_j^{-1}, & j=n.
\end{cases}
$$

The $*$-Hopf algebra $(U_q \mathfrak{g},*)$ is a quantum analog for the
universal enveloping algebra of $\mathfrak{sp}_{2n}(\mathbb R)$.

Recall a general notion \cite{Klimyk}. Let $F$ be an algebra which is also a
module over a Hopf algebra $A$. $F$ is called an $A$-module algebra if the
multiplication $m:F\otimes F\to F$ is a morphism of $A$-modules and the unit
$1\in F$ is $A$-invariant.

In \cite{9703005}, L.Vaksman with his collaborators introduced a $U_q
\mathfrak{g}$-module algebra $\mathbb C[\mathfrak{p}^-]_q$, a quantum analog
of the algebra of holomorphic polynomials on $\mathfrak{p}^-$. They followed
V.Drinfeld's approach to producing quantum analogs of function algebras by
duality. The next two theorems give an explicit description of the $U_q
\mathfrak{g}$-module algebra
 $\mathbb C[\mathfrak{p}^-]_q$ in terms of generators and relations.

%He also proved the existence and uniqueness of the faithfull
%irreducible $*$-representation of $\mathrm {Pol}(\mathfrak{p}^-)_q$
%and the $(U_q \mathfrak{g},*)$-invariant integral over the bounded
%symmetric domain $\mathbb D$ (see \cite{thesis}[chap.2]).

%We define a quantum analog of the $U \mathfrak{sp}_{2n}(\mathbb R)$-module
%$*$-algebra $\mathrm{Pol}(\mathfrak{p}^-)$ of polynomials on
%$\mathfrak{p}^-_{\mathbb R}$ and introduce a faithful irreducible
%$*$-representation of the quantum polynomial algebra. Also we investigate a
%quantum analog of a unique $U \mathfrak{sp}_{2n}(\mathbb R)$-invariant
%integral over $\mathbb D$.

%All these results are obtained by L.Vaksman with his collaborators for an
%arbitrary irreducible quantum bounded symmetric domain (see \cite{0109198}).
%A construction of the $(U_q \mathfrak{sp}_{2n},*)$-module $*$-algebra
%$\mathrm {Pol}(\mathfrak{p}^-)_q$, a quantum analog of the polynomial
%algebra on $\mathfrak{p}^-$, can be found in \cite{9703005}.
 %First of all we should define a
%$U_q \mathfrak{g}$-module algebra $\mathbb C[\mathfrak{p}^-]_q$, which is a
%q-analog of the algebra of holomorphic polynomials on $\mathfrak{p}^-$.

%Here we give an explicit description of $\mathbb C [\mathfrak{p}^-]_q$ in
%terms of generators and relations.

\begin{theorem}
$\mathbb C[\mathfrak{p}^-]_q$ is isomorphic to the algebra generated by
$z_{ij}, 1 \leq j \leq i \leq n,$ whose ddefinig relations are the
following:
\begin{align}
 &z_{ij}z_{kl}=q^2z_{kl}z_{ij}, &i=j=l < k \label{rel_1}
\\ &z_{ij}z_{kl}=q^2z_{kl}z_{ij}, &j <i=k=l
\\ &z_{ij}z_{kl}=qz_{kl}z_{ij}, &j<l<i=k
\\ &z_{ij}z_{kl}=q z_{kl}z_{ij}, &j=l<i<k
\\ &z_{kl}z_{ij}=z_{ij}z_{kl}, &j<l \leq k<i
\\ &z_{ij}z_{kl}=z_{kl}z_{ij}+q(q^2-q^{-2})z_{lj}z_{ki}, &i=j<k=l
\\ &z_{ij}z_{kl}=z_{kl}z_{ij}+(q^2-q^{-2})z_{lj}z_{ki}, &i=j<l<k
\\ &z_{ij}z_{kl}=z_{kl}z_{ij}+(q^2-q^{-2})z_{lj}z_{ki}, &j<i<k=l
\\ &z_{ij}z_{kl}=z_{kl}z_{ij}+(q-q^{-1})(qz_{li}z_{kj}+z_{ki}z_{lj}), &j<i<l<k
\\ &z_{ij}z_{kl}=z_{kl}z_{ij}+(q-q^{-1})z_{il}z_{kj}, &j<l<i<k
\\ &z_{ij}z_{kl}=qz_{kl}z_{ij}+(q-q^{-1})z_{il}z_{kj}, &j<i=l<k.
\label{rel_11}
\end{align}
\end{theorem}

\begin{theorem}
$\mathbb C[\mathfrak{p}^-]_q$ is equipped with the structure of $U_q
\mathfrak{g}$-module algebra defined on generators as follows: for $k<n$
\begin{equation*}
K_k z_{ij} =\begin{cases} q^2 z_{ij},& i=j=k,
\\ q^{-2}z_{ij}, & i=j=k+1,
\\ qz_{ij}, & i=k>j \,\text{ or }\, i-1>k=j,
\\ q^{-1}z_{ij}, & i-1=k>j \,\text{ or }\, i>k+1=j,
\\ z_{ij}, & \text{otherwise}.
\end{cases}
\end{equation*}

\begin{equation*}
E_k z_{ij}= q^{-1/2}
\begin{cases}
 (q+q^{-1})z_{i \, j-1}, & i=j=k+1,
\\ z_{i-1 \, j}, & i=k+1>j,
\\ z_{i \, j-1}, & i>k+1=j,
\\ 0, & \text{otherwise},
\end{cases}
\quad
 F_k z_{ij} = q^{1/2} \begin{cases}
(q+q^{-1})z_{i+1\,j}, \,&\, i=j=k,
\\ z_{i+1\,j}, & i=k>j,
\\ z_{i\,j+1}, & i>k=j,
\\ 0, & \text{otherwise}.
\end{cases}
\end{equation*}
and
\begin{equation*}
K_n z_{ij}=
\begin{cases}
q^4 z_{ij}, & i=j=n,
\\q^2 z_{ij}, & i=n>j,
\\ z_{ij}, & \text{otherwise}.
\end{cases}
\end{equation*}
\begin{equation*}
F_n z_{ij}= \begin{cases} q, & i=j=n,
\\ 0, & \text{otherwise}.
\end{cases}
\qquad E_n z_{ij}=-
\begin{cases}
qz_{nn}z_{ij}, & i=n \geq j,
\\ z_{ni}z_{nj}, & \text{otherwise}.
\end{cases}
\end{equation*}
\end{theorem}
% Similarly, consider the algebra $\overline{\mathbb C[\mathfrak{p}^-]_q}$ generated by
%$z^*_{ij}, 1 \leq j \leq i \leq n,$ and the corresponding relations.

Note, that the algebra $\mathbb C[\mathfrak{p}^-]_q$ was obtained in
the Kamita's paper \cite{Kamita}, but only as a $U_q
\mathfrak{k}$-module algebra.

In \cite{9703005}, a $(U_q \mathfrak{g},*)$-module $*$-algebra $\mathrm
{Pol}(\mathfrak{p}^-)_q$, (i.e. we have $ (\xi f)^*=S(\xi)^*f^*$ for all
$\xi \in U_q \mathfrak{g}, f \in \mathrm{Pol}(\mathfrak{p}^-)_q $) was
introduced. It is considered as a quantum analog of the algebra of
polynomials on $\mathfrak{p}^-_{\mathbb R}$. Also the existence and
uniqueness of the faithfull irreducible $*$-representation of $\mathrm
{Pol}(\mathfrak{p}^-)_q$ and the $(U_q \mathfrak{g},*)$-invariant integral
over the bounded symmetric domain $\mathbb D$ were proved (see
\cite[chap.2]{thesis}).

\begin{proposition}
A list of defining relations for the $*$-algebra
$\mathrm{Pol}(\mathfrak{p}^-)_q$ consists of the relations
\eqref{rel_1}-\eqref{rel_11} and
\begin{align*}
 & z_{ij}^*z_{kl}=z_{kl}z^*_{ij},
& \!\!\!\!\!\!j\neq k,l\,\, \& \,\, i\neq k,l,
\\& z_{ij}^*z_{kl}=qz_{kl}z_{ij}^*-(q^{-1}-q)\!\!\sum_{m>k}z_{ml}z_{mj}^*,
 & \!\!\!\!\!\!i=k>j>l,
\\& z_{ij}^*z_{kl}=qz_{kl}z_{ij}^*-q(q^{-1}-q)
(\!\!\sum_{i\geq m >k}z_{ml}z_{im}^* +q\!\!\sum_{m>i}z_{ml}z_{mi}^*), &
\!\!\!\!\!\!i>k=j>l,
\\& z_{ij}^*z_{kl}=qz_{kl}z_{ij}^*-(q^{-1}-q)(\!\!\sum_{k \geq m>l}z_{km}z_{im}^*
+q\!\!\sum_{i \geq m>k}z_{mk}z_{im}^* +q^2\!\sum_{m>i}z_{mk}z_{mi}^*), &
\!\!\!\!\!\!i>k>j=l,
\\& z_{ij}^*z_{kl}=q^2z_{kl}z_{ij}^*-(1+q^2)(q^{-1}-q)(\!\!\sum_{i\geq m>l}
z_{mk}z_{im}^* + q\!\!\sum_{m>i}z_{mk}z_{mi}^*), & \!\!\!\!\!\!i>j=k=l,
\\& z_{ij}^*z_{kl}=
q^2z_{kl}z_{ij}^*-(1+q^2)(q^{-1}-q)\!\!\sum_{m>k}z_{ml}z_{mi}^*, &
\!\!\!\!\!\!i=j=k>l,
\\& z_{ij}^*z_{kl}=q^2z_{kl}z_{ij}^*-q(q^{-1}-q)(\!\!\sum_{i \geq k'>j}
z_{kk'}z_{ik'}^* +\!\!\sum_{k'>i}z_{k'l}z_{k'j}^*+ q^2\!
\sum_{k'>i}z_{k'k}z_{k'i}^*)
\\& +(q^{-1}-q)^2\!\!\!\sum_{k'>i,l'>j,k'>l'}z_{k'l'}z_{k'l'}^*
+(q^{-1}-q)^2\!\sum_{k'>i}z_{k'k'}z_{k'k'}^*+1-q^2, & \!\!\!\!\!\!i=k>j=l,
\\& z_{ij}^*z_{kl}=q^4z_{kl}z_{ij}^*-q(q^{-1}-q)(1+q^2)^2 \!\sum_{k'>i}
z_{k'l}z_{k'j}^* \\& +(q^{-1}-q)^2(1+q^2)\!\!\sum_{k'>i}z_{k'k'}z^*_{k'k'}
+(q^{-1}-q)^2(1+q^2)^2 \!\!\!\sum_{k'>j'>i}z_{k'j'}z_{k'j'}^*+1-q^4, &
\!\!\!\!\!\!i=j=k=l.
\end{align*}
together with relations which can be obtained from the above due to obvious
involution properties (the involution $*$ is defined naturally $*: z_{ij}
\mapsto z^*_{ij}$).
\end{proposition}

\section{The faithfull representation and the invariant integral}

In this section we outline some results which were obtained in
\cite[chap.2]{thesis}.

Introduce an irreducible $*$-representation of
$\mathrm{Pol}(\mathfrak{p}^-)_q$. Let $\mathcal H$ be a
$\mathrm{Pol}(\mathfrak{p}^-)_q$-module with a single generator $v_0$ and
the relations
$$z^*_{ij}v_0=0, \quad 1 \leq j \leq i \leq n.$$
\begin{proposition}\cite[sect. 2.2]{thesis}
1. $\mathcal {H}=\mathbb C[\mathfrak{p}^-]_qv_0$.

2. There exists a unique sesquilinear form $(\cdot,\cdot)$ on $\mathcal H$
with the properties: i) $(v_0,v_0)=1$; ii) $(fv,w)=(v,f^*w)$ for all $v,w
\in \mathcal H$, $f \in \mathrm{Pol}(\mathfrak{p}^-)_q$.

3. The form $(\cdot,\cdot)$ is positive definite on $\mathcal H$.
\end{proposition}

Denote by $T_F$ the representation of $\mathrm{Pol}(\mathfrak{p}^-)_q$:
$$T_F(f)v=fv, \quad f \in \mathrm{Pol}(\mathfrak{p}^-)_q, v \in
\mathcal H.$$
\begin{theorem}\cite[sect. 2.2]{thesis}
1. $T_F$ is a faithful irreducible $*$-representation.

2. A $\mathrm{Pol}(\mathfrak{p}^-)_q$-representation with such properties is
unique up to unitary equivalence.
\end{theorem}

Let $d\nu$ be an invariant measure on the irreducible bounded symmetric
domain $\mathbb D$. Note that
$$
\int_{\mathbb D}fd\nu=\infty, \quad f \in \mathrm{Pol}(\mathfrak{p}^-), f
\neq 0.
$$
So we have to construct a quantum analog of the algebra of smooth functions
on $\mathbb D$ with compact supports and to define an invariant integral on
it.

The algebra $\mathbb C[\mathfrak{p}^-]_q$ is equipped with a natural grading
$$
\mathbb C[\mathfrak{p}^-]_q=\bigoplus_{k=0}^{\infty} \mathbb
C[\mathfrak{p}^-]_{q,k}, \quad \deg z_{ij}=1.
$$
%The vector space $\mathcal{H}$ inherits the grading $\mathcal{H}=\bigoplus
%\limits_{k=0}^{\infty} \mathcal{H}_k$.

Extend the $*$-algebra $\operatorname{Pol}(\mathfrak{p}^-)_q$ by attaching
an element $f_0$ which satisfies the following relations: $$f_0^2=f_0,\;
f_0^*=f_0,\quad \psi^*f_0=f_0\psi=0,\quad
 \psi\in\mathbb{C}[\mathfrak{p}^-]_{q,1}.$$

The two-sided ideal of the extended algebra $$\mathcal{D}(\mathbb
D)_q\stackrel{\rm def}{=}\mathrm{Pol}(\mathfrak{p}^-)_q \cdot f_0 \cdot
\mathrm{Pol}(\mathfrak{p}^-)_q$$ is treated as a quantum analog of the space
of smooth functions with compact supports on $\mathbb D$.

$\mathcal{D}(\mathbb D)_q$ is equipped with a $(U_q \mathfrak{g},*)$-module
algebra structure via
\begin{equation*}
F_j f_0=\begin{cases} -\frac{q^5}{1-q^4}f_0z^*_{nn},& j=n,
\\0, & j \neq n,
\end{cases}
\qquad E_j f_0=\begin{cases} -\frac{q}{1-q^4}z_{nn}f_0,& j=n,
\\0, & j \neq n.
\end{cases}
\end{equation*}

Let $\mathcal{H}_F\overset{\mathrm{def}}{=}
\mathbb{C}[\mathfrak{p}^-]_qf_0.$

$\mathcal{H}_F$ is a $\mathcal{D}(\mathbb{D})_q$-module, a
$\operatorname{Pol}(\mathfrak{p}^-)_q$-module, and a
$U_q\mathfrak{g}$-module. Denote by $\mathcal{T}_F$ the corresponding
representations of $\mathcal{D}(\mathbb{D})_q$ and
$\operatorname{Pol}(\mathfrak{p}^-)_q$ in the vector space $\mathcal{H}_F$
(these representations are related with their faithful irreducible
$*$-representations) and by $\Gamma$ the representation of $U_q
\mathfrak{g}$.

Define a linear functional on $\mathcal{D}(\mathbb D)_q$:
$$
\int_{\mathbb D_q} f
d\nu=(1-q^4)^{\frac{n(n+1)}{2}}\mathrm{tr}(\mathcal{T}_F(f) \Gamma(K^{-1})),
\quad f \in \mathcal{D}(\mathbb D)_q,
$$
with $K=K_1^{2n}K_2^{2(2n-1)}...K_{n-1}^{(n-1)(n+2)}K_n^{n(n+1)/2}$.
\begin{theorem}\cite[sect. 2.2]{thesis}
1. The above integral is $U_q \mathfrak{g}$-invariant and positive, i.e.
$$
\int_{\mathbb D_q} (\xi \cdot f) d\nu= \varepsilon(\xi) \int_{\mathbb D_q} f
d\nu, \quad \xi \in U_q \mathfrak{g}, f \in \mathcal{D}(\mathbb D)_q,
$$
and
$$
\int_{\mathbb D_q} (f^*f) d\nu>0, \quad f \in \mathcal{D}(\mathbb D)_q, f
\neq 0.
$$
2. A positive $U_q \mathfrak{g}$-invariant integral on $\mathcal{D}(\mathbb
D)_q$ is unique up to a constant multiple.
\end{theorem}

\section*{Acknowledgment}

The author would like to thank professor L.Vaksman for statements of the
problems and helpful discussions.

\end{document}